\documentclass[12pt]{amsart}

\usepackage{amsmath,xspace,amssymb,mathrsfs}
\usepackage{color}

\input xy
\xyoption{all}
\xyoption{2cell}
\UseAllTwocells
\CompileMatrices

\newcommand{\Spec}{\operatorname{Spec}}
\renewcommand{\phi}{\varphi}

\newcommand{\Ker}{\operatorname{Ker}}
\newcommand{\Ima}{\operatorname{Im}}

\newcommand{\V}{\operatorname{V}}

\newcommand{\Sp}{\operatorname{Sp}}

\newcommand{\colimit}{\operatorname{colim}}

\newtheorem{proposition}{Proposition}[section]
\newtheorem{lemma}[proposition]{Lemma} 
\newtheorem{corollary}[proposition]{Corollary}
\newtheorem{theorem}[proposition]{Theorem}

\newtheorem{prop-def}[proposition]{Proposition and definition}

  \theoremstyle{definition}
  \newtheorem{definition}[proposition]{Definition}

  \newtheorem{remark}[proposition]{Remark}

\begin{document}

\title{Flat topology and its dual aspects}

\author{Abolfazl Tarizadeh}
\address{ Department of Mathematics, Faculty of Basic Sciences, University of Maragheh \\
P. O. Box 55136-553, Maragheh, Iran.
 }
\email{ebulfez1978@gmail.com}


\begin{abstract} In this article, a new and natural topology on the prime spectrum is introduced which behaves completely as the dual of the Zariski topology. It is called the flat topology. The basic and also some sophisticated properties of the flat topology are proved. Specially, various algebraic characterizations for the noetherianness of the flat topology is given. Using the flat topology, then some facts on the structure of the prime ideals of a ring come to light which are not in the access of the Zariski topology.
\end{abstract}

\maketitle

\section{Introduction}

The main purpose of the present article is to go to the beyond of the Zariski topology by establishing a new topology on the prime spectrum in order to bring under the light some facts on the structure of the prime ideals of a ring which are not possible to discover them by using the Zariski topology. We may analogize the Zariski topology as one of the wings of a bird. A bird in order to fly needs the two wings. The new topological structure on the prime spectrum, which we are going to establish it and will be act as the second wing of a bird, will be called the \emph{flat topology}.  Our starting point for establishing the flat topology originates from a simple observation. In fact it emanates from the structure of the constructible topology. The flat topology, as we shall observe it in the article, behaves completely as the dual of the Zariski topology even though the initial construction of the flat topology seems far from to be the dual of the Zariski topology. Because, in the flat topology, the closed subsets of $\Spec(R)$ are coming from the flat $R-$algebras, see Proposition \ref{prop- 55}. But, by applying the going-down property of the flat algebras, it is then proved that the flat topology is the coarsest topology on the prime spectrum $\Spec(R)$ such that the subsets $D(f)=\{\mathfrak{p}\in\Spec(R) : f\notin\mathfrak{p}\}$ are closed where $f$ runs through the underlying set of $R$, see Theorem \ref{remark 122}. While the Zariski topology, as we know, is the coarsest topology such that the subsets $D(f)$ are open. This ensures that the flat topology could be a suitable candidate for being the dual of the Zariski topology. The opens of the flat topology, unlike the Zariski, are sufficiently small. But, as the Zariski, the flat topology is also spectral.\\

``Localization" and ``quotient'' are two fundamental and vital tools in commutative algebra. The Zariski closed subsets has the ``quotient" nature while the flat closed subsets has the ``localization" nature. For example, if $\mathfrak{p}$ is a prime ideal of a ring $R$ then the Zariski closure of the point $\mathfrak{p}$ in $\Spec(R)$ comes from the canonical ring map $R\rightarrow R/\mathfrak{p}$ while the flat closure of this point comes from the canonical ring map $R\rightarrow R_{\mathfrak{p}}$, see Corollary \ref{coro 1}. In fact, there is a complete duality in the topological level between the localizations and quotients; as an evidence for this see Theorem \ref{coro 12453}.\\

We should mention that the flat topology and  Hochster's inverse topology \cite[Proposition 8]{Hochster} are exactly the same things, see Remark \ref{Remark 0020}. But our approach to construct the flat topology is completely different from his method. Our approach only uses the standard results and methods of commutative algebra. In subsequent works we give some applications of the flat topology. Specially, it is proved that the set of minimal primes of a ring is quasi-compact w.r.t. the flat topology. This result easily implies some of the major results in the literature on the finiteness of minimal primes. There is also another remarkable application of the flat topology which states that ``the finitely generated flat modules of a ring $R$ are projective if and only if every flat closed subset of $\Spec(R)$ which is stable under the specialization is flat open'', see \cite[Theorem 3.3]{Tarizadeh 2020}. This result, in particular, implies that if a ring $R$ has either a finitely many minimal primes or a finitely many maximal ideals then every finitely generated  flat $R-$module is projective, see \cite[Theorem 3.4]{Tarizadeh 2020}. The latter result, as stated in \cite{Tarizadeh 2020}, vastly generalizes some major results in the literature on the projectivity of finitely generated flat modules. \\

After that a primary version of the present article was submitted to the arXiv preprint server, the author received an e-mail from Professor Marco Fontana. In that letter, in addition to expressing of his interest to our work,
he also informed us about their earlier work \cite{Fontana}. During the writing the present article the author was not aware of the existence of such article in the literature. The author then went through that article and found some of the overlappings (surprisingly there is also an overlapping on the terminology of  the ``\emph{flat topology}") between their article and our work which had been obtained independently. Their article, after Fontana's informing, was so useful to the author in the understanding the structure of the flat topology more deeply. But there is a gap in that article where they deduce that the connected components of the prime spectrum behave quite differently in the passage from the Zariski topology to the flat topology \cite[Remark 2.20, (f)]{Fontana}. This is a wrong statement. By the way, the flat and Zariski connected components of the prime spectrum are precisely the same, see Theorem \ref{prop 9}. \\

Note that in spite of the existence of some similarities between the flat and Zariski topologies such as their connected components and Hausdorfness, there are in general tremendous differences between the two topologies. In fact, the two topologies on $\Spec(R)$ are the same if and only if every prime ideal of $R$ is maximal \cite[Proposition 2.3]{Ebulfez 25}. The noetheriannes of the infinite spectra is quite different in the passage between the two topologies, see Theorem \ref{theorem 4537}. \\

In summary, the readers will find the flat topology very interesting and efficient for many applications. Throughout the article, all of the rings are commutative. \\

The term ``flat topology" has different meanings in the literature (where it is sometimes also called the fppf topology or fpqc as well; these are Grothendieck topologies). But we point out that, in this article, our purpose of the ``flat topology" is actually different from these concepts. \\

The titles of the sections should be sufficiently explanatory, but we explain them a little further. In Section 2, the constructible topology is introduced and its initial properties are studied. In Section 3, first the flat topology is established and then its basic properties are extracted. Specially, the going-down property of the flat algebras is used to produce a basis for the opens of the flat topology (Theorem \ref{remark 122}). It is proved that the flat topology is an ideal-theoretic structure (Theorem \ref{prop 10}). The irreducible and connected components of the flat topology are characterized (Corollary \ref{prop 8} and Theorem \ref{prop 9}). In the final section, various algebraic characterizations for the noetherianness of the flat topology are given, see Theorem \ref{theorem 1}. Chevalley's Theorem \cite[Tag 00FE ]{Johan}  plays a major role in the backstage of this characterization.\\

\section{The constructible topology}

The patch topology was introduced by Hochster \cite[\S2]{Hochster} for general topological spaces. In this section we re-cover the patch topology on the prime spectrum by using purely algebraic methods. Our approach is completely different from the Hochster's method. \\

Let $\phi:R\rightarrow A$ be a ring map and $\mathfrak{p}$ a prime ideal of $R$. Then $\mathfrak{p}$ is in the image of the induced map $\phi^{\ast}:\Spec(A)\rightarrow\Spec(R)$ if and only if $A\otimes_{R}\kappa(\mathfrak{p})\neq0$ where $\kappa(\mathfrak{p})$ is the residue field of $R$ at $\mathfrak{p}$. Because if $A\otimes_{R}\kappa(\mathfrak{p})\neq0$ then choose a prime ideal $P$ of it.
Let $\mathfrak{q}:=\lambda^{-1}(P)$ where $\lambda:A\rightarrow A\otimes_{R}\kappa(\mathfrak{p})$ is the canonical ring map. Then $\phi^{-1}(\mathfrak{q})=(\mu\circ\pi)^{-1}(P)=\mathfrak{p}$ where  $\mu:\kappa(\mathfrak{p})\rightarrow A\otimes_{R}\kappa(\mathfrak{p})$ and $\pi:R\rightarrow\kappa(\mathfrak{p})$ are the canonical ring maps. Conversely, if $\mathfrak{p}\in\Ima\phi^{\ast}$ there there exists a prime ideal $\mathfrak{q}$ of $A$ such that $\phi^{-1}(\mathfrak{q})=\mathfrak{p}$. There exists a ring map $\theta:A\otimes_{R}\kappa(\mathfrak{p})
\rightarrow\kappa(\mathfrak{q})$ such that $\theta\circ\lambda:A\rightarrow\kappa(\mathfrak{q})$ is the canonical ring map where $\kappa(\mathfrak{q})$ is the residue field of $A$ at $\mathfrak{q}$.
It follows that $A\otimes_{R}\kappa(\mathfrak{p})\neq0$ because $\kappa(\mathfrak{q})\neq0$.\\

\begin{lemma}\label{tensor lemma} Let $A_{1},...,A_{p}$ be $R-$algebras with the structure morphisms $\phi_{i}: R\rightarrow A_{i}$. Then $\Ima\phi^{\ast}=
\bigcap\limits_{i=1}^{p}\Ima\phi_{i}^{\ast}$ where $\phi: R\rightarrow A_{1}\otimes...\otimes A_{p}$ is the push-out of the $\phi_{i}$.  \\
\end{lemma}

{\bf Proof.} Let $A:=A_{1}\otimes_{R}...\otimes_{R}A_{p}$. The inclusion
$\Ima\phi^{\ast}\subseteq
\bigcap\limits_{i=1}^{p}\Ima\phi^{\ast}_{i}$ is obvious. Conversely, suppose $\mathfrak{p}\in\bigcap\limits_{i=1}^{p}\Ima\phi_{i}^{\ast}$. Then $A_{i}\otimes_{R}\kappa(\mathfrak{p})$ is a non-zero vector space over $\kappa(\mathfrak{p})$ for all $i$. Therefore
$$\big(A_{1}\otimes_{R}\kappa(\mathfrak{p})\big)
\otimes_{\kappa(\mathfrak{p})}...\otimes_{\kappa(\mathfrak{p})}
\big(A_{p}\otimes_{R}\kappa(\mathfrak{p})\big)\neq0$$ (recall that if $V$ and $W$ are two non-zero vector spaces over a field $K$ then $V\otimes_{K}W\neq0$). This implies that $A\otimes_{R}\kappa(\mathfrak{p})\neq0$. Thus
$\mathfrak{p}\in\Ima\phi^{\ast}$.   $\Box$ \\

\begin{proposition}\label{direct limit} Let $(A_{i}, \psi_{ij})$ be an inductive (direct) system of $R-$algebras with the structure morphisms $\phi_{i}: R\rightarrow A_{i}$ over the directed set $I$. Then $\Ima\phi^{\ast}=\bigcap\limits_{i\in I}\Ima\phi_{i}^{\ast}$
where $\phi: R\rightarrow\colimit\limits_{i\in I}A_{i}$ is induced by the $\phi_{i}$.\\
\end{proposition}

{\bf Proof.} The inclusion $\Ima\phi^{\ast}\subseteq\bigcap\limits_{i\in I}\Ima\phi_{i}^{\ast}$  is obvious. For the reverse inclusion, pick  $\mathfrak{p}\in\bigcap\limits_{i\in I}\Ima\phi_{i}^{\ast}$. It suffices to show that $\big(\colimit\limits_{i\in I}A_{i}\big)\otimes_{R}\kappa(\mathfrak{p})\neq0$. We have the canonical isomorphism $\big(\colimit\limits_{i\in I}A_{i}\big)\otimes_{R}\kappa(\mathfrak{p})
\simeq\colimit\limits_{i\in I}\big(A_{i}\otimes_{R}\kappa(\mathfrak{p})\big)$.
But $A_{i}\otimes_{R}\kappa(\mathfrak{p})\neq0$ for all $i\in I$. Therefore  their inductive limit is also nontrivial. $\Box$  \\

\begin{proposition}\label{constructible} Let $R$ be a ring. Then the collection of subsets $\Ima\phi^{\ast}$ of $\Spec(R)$ where $\phi: R\rightarrow A$ is a ring homomorphism satisfies in the axioms of closed subsets in a topological space. \\
\end{proposition}

{\bf Proof.} Let $\phi_{i}:R\rightarrow A_{i}$ be a finite family of ring maps with $i=1,...,n$. Consider the ring map $\phi:R\rightarrow A:=A_{1}\times...\times A_{n}$ which maps each $r\in R$ into $\big(\phi_{1}(r),...,\phi_{n}(r)\big)$. Then $\Ima\phi^{\ast}=\bigcup\limits_{i=1}^{n}\Ima\phi^{\ast}_{i}$. Because, every prime ideal of $A$ is of the form $A_{1}\times...\times A_{j-1}\times\mathfrak{p}_{j}\times A_{j+1}\times...\times A_{n}$ where $\mathfrak{p}_{j}$ is a prime ideal of $A_{j}$ for some $j$. Therefore the collection is stable under the finite unions. Now let $\{\Ima\phi_{i}^{\ast}\}_{i\in I}$ be a subset of the collection where for each $i$, $\phi_{i}: R\rightarrow A_{i}$ is a ring map. For each finite subset $J$ of  $I$, let $A_{J}=\bigotimes\limits_{i\in J}A_{i}$ be the tensor product of the $R-$algebras $A_{i}$ with $i\in J$ and let $\phi_{_{J}}: R\rightarrow A_{J}$ be the structure morphism of the $R-$algebra $A_{J}$ induced by the $\phi_{i}$.
Note that if $J$ is the empty set then $A_{J}=R$ and if $J=\{i\}$ for some $i$ then $A_{J}=A_{i}$. For each two finite subsets $J$ and $J'$ of $I$ with $J\subseteq J'$, let $\psi_{_{JJ'}}: A_{J}\rightarrow A_{J'}$ be the canonical ring map which maps each pure tensor $\bigotimes\limits_{i\in J}a_{i}$ of $A_{J}$ into $\bigotimes\limits_{i\in J'}b_{i}$ where for each $i\in J$, $b_{i}=a_{i}$ and for each $i\in J'\setminus J$, $b_{i}=e_{i}$ is the identity element of $A_{i}$. Therefore $(A_{J}, \psi_{_{JJ'}})$ forms an inductive system of $R-$algebras over the set of finite subsets of $I$ which is a directed poset ordered by the inclusion. Let
$\phi: R\rightarrow\colimit\limits_{J}A_{J}$ be the structure morphism induced by the $\phi_{_{J}}$.
Then by Proposition \ref{direct limit}, $\Ima\phi^{\ast}=
\bigcap\limits_{J}\Ima\phi_{_{J}}^{\ast}$.
Also, by Lemma \ref{tensor lemma}, $\Ima\phi_{_{J}}^{\ast}=\bigcap\limits_{i\in J}\Ima\phi_{i}^{\ast}$. Therefore $\Ima\phi^{\ast}=\bigcap\limits_{i\in I}\Ima\phi^{\ast}_{i}$. Thus the collection is also stable under the arbitrary intersections. $\Box$ \\

The resulted topology in Proposition \ref{constructible} is called the \emph{constructible} topology.\\

Recall that a topological space is said to be compact if it is quasi-compact and Hausdorff. Also recall that a topological space is said to be totally disconnected if its connected subsets are just the single-point subsets. Every compact and totally disconnected topological space is called a profinite space.\\

Note that if we consider the patch topology \cite[\S2]{Hochster} on the prime spectrum $\Spec(R)$ then the collection of subsets $D(f)\cap V(I)$ where $f\in R$ and $I$ is a finitely generated ideal of $R$ formes a basis for the opens of the patch topology. It follows that the patch topology is coarser than of the constructible topology for now. Therefore we are not able to apply the compactness of the patch topology  \cite[Theorem 1]{Hochster} to deduce the following result.  It requires a direct proof: \\

\begin{theorem}\label{prop 1} The set $\Spec(R)$ equipped with the constructible topology is a profinite space.\\
\end{theorem}

{\bf Proof.} The constructible topology is clearly Hausdorff. To prove its quasi-compactness, let $\mathscr{C}=\{\Ima\phi_{i}^{\ast}\}_{i\in I}$ be a family of closed subsets of  $\Spec(R)$ with the finite intersection property where $\phi_{i}: R\rightarrow A_{i}$ is a ring homomorphism for all $i$. With taking into account the notations which used in the proof of Proposition \ref{constructible}, then for each finite subset $J$ of $I$, the $R-$algebra $A_{J}$ is a nontrivial ring since the family $\mathscr{C}$ has the finite intersection property. Therefore  $\colimit\limits_{J}A_{J}$ is a nontrivial ring. This, in particular, implies that $\bigcap\limits_{i\in I}\Ima\phi^{\ast}_{i}\neq\emptyset$. Finally, let $C$ be a connected component of $\Spec(R)$. There exists an $R-$algebra $A$ with the structure morphism $\phi: R\rightarrow A$ such that $C=\Ima\phi^{\ast}$. Let $\mathfrak{p}$ and $\mathfrak{p'}$ be two prime ideals of $A$. It suffices to show that $\phi^{\ast}(\mathfrak{p})=\phi^{\ast}(\mathfrak{p'})$. Suppose on the contrary, then we may choose an element $f\in\phi^{\ast}(\mathfrak{p})$ such that $f\notin\phi^{\ast}(\mathfrak{p'})$. This implies that $\phi^{\ast}(\mathfrak{p})\in V(f)\cap C$ and $\phi^{\ast}(\mathfrak{p'})\in D(f)\cap C$. Therefore $V(f)\cap C$ and $ D(f)\cap C$ are nonempty. They are also disjoint and open subsets of $C$. But this is a contradiction since $C$ is connected. $\Box$ \\

\begin{remark}\label{remark 432} The constructible topology is the coarsest topology over $X=\Spec(R)$ such that for each $f\in R$ the subset $D(f)$ is both open and closed. Because, let $\mathscr{T}$ be the such topology and let $\mathscr{T}'$ be the constructible topology. Clearly $\mathscr{T}\subseteq\mathscr{T}'$. So the identity map $Id: (X,\mathscr{T}')\rightarrow(X,\mathscr{T})$ is continuous. It is also a closed map since $\mathscr{T}$ is Hausdorff and $\mathscr{T}'$ is quasi-compact, see Theorem \ref{prop 1}. Therefore $Id$ is a homeomorphism and so $\mathscr{T}=\mathscr{T}'$. It follows that the collection of subsets $D(f)\cap V(g)$ with $f,g\in R$ formes a sub-basis for the opens of the constructible topology. In other words, the collection of subsets $D(f)\cap\V(I)$ where $f\in R$ and $I$ is a finitely generated ideal of $R$ is a basis for the opens of the constructible topology. Therefore the constructible and patch topologies of $\Spec(R)$ are the same. \\
\end{remark}

\section{The flat topology}

In this section, completely inspired from the structure of the constructible topology, we establish a new topology on the prime spectrum which will be called the flat topology. Then its basic properties are studied.\\

\begin{proposition}\label{prop- 55} Let $R$ be a ring. Then the collection of subsets $\Ima\phi^{\ast}$ of $\Spec(R)$ where $\phi: R\rightarrow A$ is a flat ring homomorphism satisfies in the axioms of closed subsets in a topological space.\\
\end{proposition}

{\bf Proof.} This is proven exactly like Proposition \ref{constructible} with taking into account the facts that the tensor product of flat modules and also the inductive (direct) limit of every inductive system of flat modules are flat. $\Box$ \\

The resulted topology of Proposition \ref{prop- 55} is called the flat topology. By comparing Proposition \ref{prop- 55} with Proposition \ref{constructible} the initial difference between the flat and constructible topologies is emerged, i.e., the flat topology is coarser than of the constructible topology.\\

A subset $E$ of $\Spec(R)$ is said to be a flat closed (resp. Zariski closed, constructible closed) if it is a closed subset of $\Spec(R)$ with respect to the flat (resp. Zariski, constructible) topology. Throughout the article, we shall use freely these phrases even more phrases such as flat open, flat irreducible, Zariski open and etc with the appropriate meanings. \\

Note that if we consider the inverse topology \cite[Prop. 8]{Hochster} on the prime spectrum $\Spec(R)$ then the collection of subsets $V(f)$ with $f\in R$ formes a sub-basis for the opens of the inverse topology. It follows that the inverse topology is coarser than of the flat topology for now. The following result says that these two topologies are the same things. \\

\begin{theorem}\label{remark 122} The collection of subsets $V(I)$ where $I$ runs through the set of finitely generated ideals of $R$ forms a basis for the flat opens of $\Spec R$. \\
\end{theorem}

{\bf Proof.} For each $f\in R$ then the canonical map $\pi:R\rightarrow R_{f}$ is a flat ring map. Therefore $D(f)=\Ima\pi^{\ast}$ is a flat closed subset of $\Spec R$ and so $V(f)=\Spec R\setminus D(f)$ is a flat open subset of $\Spec(R)$. It follows that if $I=(f_{1},...,f_{n})$ is a finitely generated ideal of $R$ then $V(I)=\bigcap\limits_{i=1}^{n}V(f_{i})$ is a flat open subset of $\Spec R$. Conversely, let $U$ be an arbitrary flat open subset of $\Spec R$. Then there exists a flat ring map $\phi:R\rightarrow A$ such that $U=\Spec R\setminus\Ima\phi^{\ast}$. Thus, if $\mathfrak{p}$ is a prime ideal of $R$ then $\mathfrak{p}\in U$ if and only if $\mathfrak{p}A=A$. Because every flat ring map has the going-down property. Therefore, if $\mathfrak{p}\in U$ then there are a finitely many elements $g_{1},...,g_{n}\in\mathfrak{p}$ such that $1_{A}=\sum\limits_{i=1}^{n}\phi(g_{i})a_{i}$ where $a_{i}\in A$ for all $i$.
Let $J=(g_{1},...,g_{n})$.
Then clearly $\mathfrak{p}\in V(J)\subseteq U$. $\Box$ \\

\begin{remark}\label{Remark 0020} By Theorem \ref{remark 122}, the flat and inverse topologies are the same things. \\
\end{remark}

\begin{proposition}\label{prop 2} Let $\psi: R\rightarrow R'$ be a ring map. Then the induced map $\psi^{\ast}:\Spec(R')\rightarrow\Spec(R)$ is continuous with respect to the flat (resp. constructible) topology. \\
\end{proposition}

{\bf Proof.} For each $f\in R$ and for every finitely generated ideal $I$ of $R$, we have $(\psi^{\ast})^{-1}\big(D(f)\big)=D\big(\psi(f)\big)$ and $(\psi^{\ast})^{-1}\big(\V(I)\big)=\V(I^{e})$ where $I^{e}$ is the extension of $I$ under $\psi$ which is a finitely generated ideal of $R'$. $\Box$ \\

\begin{remark}\label{remark 12} The set $\Spec(R)$ equipped with the flat topology is quasi-compact since the flat topology is coarser than of the constructible topology. Then apply Theorem \ref{prop 1}. \\
\end{remark}

Let $E$ be a subset of $\Spec(R)$. The closure of $E$ in $\Spec(R)$ with respect to the flat topology is denoted by $\Lambda(E)$. If $E=\{\mathfrak{p}\}$ for some prime ideal $\mathfrak{p}$ of $R$ then its flat closure is simply denoted by $\Lambda(\mathfrak{p})$. \\

\begin{corollary}\label{coro 1} Consider the flat topology over $\Spec(R)$. Then for each prime ideal $\mathfrak{p}$ of $R$, $\Lambda(\mathfrak{p})=\{\mathfrak{p'}\in\Spec(R) : \mathfrak{p'}\subseteq\mathfrak{p}\}$. In particular, $\mathfrak{p}$ is a closed point of $\Spec(R)$ if and only if it is a minimal prime ideal of $R$.   \\
\end{corollary}

{\bf Proof.} It implies from Theorem \ref{remark 122}. $\Box$ \\

\begin{lemma}\label{key lemma} Let $K$ be a quasi-compact subset of $\Spec(R)$ with respect to the Zariski topology. Then $\Lambda(K)=\bigcup\limits_{\mathfrak{p}\in K}\Lambda(\mathfrak{p})$.\\
\end{lemma}

{\bf Proof.} The inclusion $\bigcup\limits_{\mathfrak{p}\in K}\Lambda(\mathfrak{p})\subseteq\Lambda(K)$ holds in every topological space. Suppose there is some $\mathfrak{p'}\in\Lambda(K)\setminus\bigcup\limits_{\mathfrak{p}\in K}\Lambda(\mathfrak{p})$. Then by Corollary \ref{coro 1}, corresponding with each $\mathfrak{p}\in K$ there exists an element $c_{\mathfrak{p}}\in\mathfrak{p'}$ such that $c_{\mathfrak{p}}\notin\mathfrak{p}$.  Thus $K\subseteq\bigcup\limits_{\mathfrak{p}\in K}D(c_{\mathfrak{p}})$. But $K$ is quasi-compact therefore we find a finite set $\{D(c_{\mathfrak{p}_{i}}): 1\leq i\leq n\}$ of elements of the covering which again covers $K$. The subset $\V(c_{\mathfrak{p}_{1}},...,c_{\mathfrak{p}_{n}})$ is a flat open neighborhood of $\mathfrak{p'}$. Therefore $\V(c_{\mathfrak{p}_{1}},...,c_{\mathfrak{p}_{n}})\cap K\neq\emptyset$. But this is a contradiction since $K\subseteq\bigcup\limits_{i=1}^{n}D(c_{\mathfrak{p}_{i}})$. $\Box$ \\

An analogue of Lemma \ref{key lemma} (even with a similar proof) holds for the Zariski closure:\\

\begin{lemma} Let $K$ be a quasi-compact subset of $\Spec(R)$ with respect to the flat topology. Then $\overline{K}=\bigcup\limits_{\mathfrak{p}\in K}\V(\mathfrak{p})$. $\Box$ \\
\end{lemma}

\begin{lemma}\label{lem 309} Let $\phi:R\rightarrow A$ be an injective ring map and let $\mathfrak{p}$ be a minimal prime ideal of $R$. Then $A$ has a (minimal) prime ideal lying over $\mathfrak{p}$. \\
\end{lemma}

{\bf Proof.} By the universal property of the localization there is a (unique) ring map $\phi_{\mathfrak{p}}:R_{\mathfrak{p}}\rightarrow S^{-1}A$ given by $r/s\rightsquigarrow\phi(r)/\phi(s)$ such that the following diagram is commutative $$\xymatrix{
R \ar[r]^{\phi} \ar[d]^{} & A \ar[d]^{} \\ R_{\mathfrak{p}}\ar[r]^{\phi_{\mathfrak{p}}} & S^{-1}A} $$ where $S=\phi(R\setminus\mathfrak{p})$.
The ring $S^{-1}A$ is non-trivial since $\phi_{\mathfrak{p}}$ is injective. Hence, $S^{-1}A$ has at least a prime ideal $S^{-1}\mathfrak{q}$. By the minimality of $\mathfrak{p}$ we have $\mathfrak{p}R_{\mathfrak{p}}=
\phi^{\ast}_{\mathfrak{p}}(S^{-1}\mathfrak{q})$. Finally, by the commutativity of the above diagram we have $\mathfrak{p}=\phi^{\ast}(\mathfrak{q})$. $\Box$ \\

\begin{definition} A subset $E$ of $\Spec(R)$ is said to be stable under the \emph{generalization} (resp. \emph{specialization}) if  for any two prime ideals $\mathfrak{p}$ and $\mathfrak{q}$ of $R$ with $\mathfrak{p}\subset\mathfrak{q}$ (resp. $\mathfrak{q}\subset\mathfrak{p}$) and $\mathfrak{q}\in E$, then $\mathfrak{p}\in E$.\\
\end{definition}

\begin{theorem}\label{prop 10} Let $E$ be a subset of $\Spec(R)$. Then the following conditions hold.\\
$\mathbf{(i)}$ The subset $E$ is flat closed if and only if it is constructible closed and stable under the generalization.\\
$\mathbf{(ii)}$ The subset $E$ is Zariski closed if and only if it is constructible closed and stable under the specialization.\\
\end{theorem}

{\bf Proof.}  $\mathbf{(i)}$: By the going-down property of the flat algebras, the implication ``$\Rightarrow$" is clear. Conversely, we have $E=\bigcup\limits_{\mathfrak{p}\in E}\Lambda(\mathfrak{p})$ since $E$ is stable under the generalization. The set $E$ is quasi-compact with respect to the Zariski topology. Therefore, by Lemma \ref{key lemma}, it is flat closed.\\
$\mathbf{(ii)}$: The implication ``$\Rightarrow$" is obvious. For the reverse implication, suppose $E=\Ima(\phi^{\ast})$ is stable under the specialization where $\phi:R\rightarrow A$ is a ring map. Let $\overline{\phi}:R/I\rightarrow A$ be the injective ring map induced by $\phi$ where $I=\Ker(\phi)$. We claim that $E=\V(I)$. The inclusion $E\subseteq\V(I)$ is obvious. Conversely, pick  $\mathfrak{p}\in\V(I)$ and let $\mathfrak{q}$ be a minimal prime ideal of $I$ such that $\mathfrak{q}\subseteq\mathfrak{p}$.
By Lemma \ref{lem 309}, there is a prime ideal $\mathfrak{q'}$ in $A$ such that $\mathfrak{q}=\phi^{\ast}(\mathfrak{q'})$. But $\mathfrak{p}\in E$ since it is stable under specialization. $\Box$ \\

In a topological space a subset which is both closed and open is called a clopen.\\

\begin{corollary} A subset of $\Spec(R)$ is flat clopen if and only if it is Zariski clopen. Moreover, the map $e\rightsquigarrow\V(e)$ is a bijection between the set of idempotent elements of $R$ and the set of flat clopens of $\Spec(R)$.  \\
\end{corollary}

{\bf Proof.} The first assertion implies from Theorem \ref{prop 10}. For the second assertion see \cite[Tag 00EE]{Johan}. $\Box$ \\

\begin{corollary}\label{coro 2} Consider the flat topology over $\Spec(R)$. Then $\Spec(R)$ is connected if and only if $R$ has no nontrivial idempotents. $\Box$ \\
\end{corollary}

The following lemma, in particular, paves the way in order to characterize the flat irreducible components of the prime spectrum.\\

\begin{lemma}\label{lemma 4465} Every irreducible and closed subset of  $\Spec(R)$ w.r.t. the flat topology has a unique generic point.\\
\end{lemma}

{\bf Proof.} Let $Z$ be an irreducible and closed subset of $\Spec(R)$ w.r.t. the flat topology. Let $\mathfrak{B}$ be the collection of non-empty subsets of the form $Z\cap\V(I)$ where $I$ is a finitely generated ideal of $R$. By Theorem \ref{remark 122}, the collection $\mathfrak{B}$ is a basis for the open subsets of $Z$ where the topology over $Z$ is induced by the flat topology. Moreover the basis $\mathfrak{B}$ has the finite intersection property since $Z$ is irreducible. We have $\bigcap\limits_{B\in\mathfrak{B}}B\neq\emptyset$. Because, suppose on the contrary, then $Z$ can be written as the union of the subsets $Z\cap\big(\Spec(R)\setminus\V(I)\big)$
where $I$ runs through the set of finitely generated ideals of $R$ such that $Z\cap V(I)\neq\emptyset$. Note that $Z$ is constructible closed and so it is quasi-compact. Moreover, for each ideal $I$ of $R$,
$\Spec(R)\setminus\V(I)$ is constructible open. Therefore there is a finite number of finitely generated ideals $I_{j}$ of $R$ with $1\leq j\leq n$, such that $Z\cap\V(I_{j})\neq\emptyset$ and
$Z\subseteq\bigcup\limits_{j=1}^{n}\Spec(R)\setminus\V(I_{j})$.
This implies that $\bigcap\limits_{j=1}^{n}Z\cap V(I_{j})=\emptyset$. But this is in contradiction with the finite intersection property of the family $\mathfrak{B}$. Therefore $\bigcap\limits_{B\in\mathfrak{B}}B\neq\emptyset$. Choose $\mathfrak{p}\in\bigcap\limits_{B\in\mathfrak{B}}B$, then we have $\Lambda(\mathfrak{p})=Z$. The uniqueness of the generic point follows from Corollary \ref{coro 1}. $\Box$ \\

The above lemma has the following immediate consequence:\\

\begin{corollary}\label{prop 8} Consider the flat topology over $\Spec(R)$. Then the map $\mathfrak{p}\rightsquigarrow\Lambda(\mathfrak{p})$ is a bijection between $\Spec(R)$ and the set of irreducible and closed subsets of  $\Spec(R)$. Moreover, under this correspondence, the set of maximal ideals of $R$ is in bijection with the set of flat irreducible components of $\Spec(R)$.  $\Box$ \\
\end{corollary}

\begin{remark} The flat irreducible components of $\Spec(R)$, in general, are different from the Zariski irreducible components. As a specific example, $\Spec(\mathbb{Z})$ is irreducible w.r.t. the Zariski topology while for each prime number $p$, the set $\big\{\{0\}, p\mathbb{Z}\big\}$ is an irreducible component of $\Spec(\mathbb{Z})$ w.r.t. the flat topology. In fact in $\Spec(R)$, the set of flat irreducible components is equal to the set of Zariski irreducible components if and only if every prime ideal of $R$ is contained in a unique maximal ideal and contains a unique minimal prime. In spite of the difference between the irreducible components but their connected components are exactly the same:\\
\end{remark}

\begin{theorem}\label{prop 9} Let  $C$ be a subset of $\Spec(R)$. Then the following conditions are equivalent.\\
$\mathbf{(a)}$ The set $C$ is a connected component of $\Spec(R)$ w.r.t. the flat topology.\\
$\mathbf{(b)}$ The set $C$ is a connected component of $\Spec(R)$ w.r.t. the Zariski topology.\\
$\mathbf{(c)}$ There exists a unique max-regular ideal $J$ of $R$ such that $C=\V(J)$.\\
\end{theorem}

To prove the above theorem we need Lemmata \ref{pierce lemma} and \ref{goodlem}. \\

Here an ideal of $R$ is said to be a regular ideal if it is generated by a subset of idempotent elements of $R$. The set of idempotent elements of $R$ is denoted by $I(R)$. Each maximal element of the set of proper regular ideals of $R$ (ordered by inclusion) is called a max-regular ideal of $R$. By the Zorn's Lemma, every proper regular ideal of $R$ is contained in a max-regular ideal of $R$. The set of max-regular ideals of $R$ is denoted by $\Sp(R)$ and it is called the \emph{pierce spectrum} of $R$. There is a profinite topology on the pierce spectrum:\\

\begin{lemma}\label{pierce lemma} The set $\Sp(R)$ is a profinite space.\\
\end{lemma}

{\bf Proof.} There is a (unique) topology over $\Sp(R)$ such that the collection of subsets $U_{e}=\{J\in\Sp(R) : e\notin J\}$ where $e$ runs through the set $I(R)$ is a basis for the opens. Because $U_{1}=\Sp(R)$ and $U_{e}\cap U_{e'}=U_{ee'}$ for all $e,e'\in I(R)$. We show that it is a profinite topology. If $J, J'\in\Sp(R)$ with $J\neq J'$ then there is an idempotent element
$e\in J$ such that $e\notin J'$. Therefore $J'\in U_{e}$ and $J\in U_{1-e}$. Note that $1-e$ is an idempotent element and $U_{1-e}=\Sp(R)\setminus U_{e}$.
This implies at once that $\Sp(R)$ equipped with this topology is Hausdorff and totally disconnected.  Now consider the map
$\psi: \Spec(R)\rightarrow \Sp(R)$ which is defined as $\psi(\mathfrak{p})=\big\langle e : e\in \mathfrak{p}\cap I(R)\big\rangle$. It is easy to check that for each prime ideal $\mathfrak{p}$, $\psi(\mathfrak{p})$ is a max-regular ideal of $R$ and so $\psi$ is well-defined. The map $\psi$ is onto.
Because let $J$ be a max-regular ideal of $R$. There is a prime ideal $\mathfrak{p}$ of $R$ such that $J\subseteq\mathfrak{p}$ since $J\neq R$. If $e\in\mathfrak{p}$ is an idempotent element then the regular ideal $J+Re$ is proper and so $e\in J$. It follows that $\psi(\mathfrak{p})=J$. The map $\psi$ is continuous w.r.t. the both Zariski and flat topologies since  $\psi^{-1}(U_{e})=D(e)=\V(1-e)$ for all $e\in I(R)$. This, in particular, implies that $\Sp(R)$ is quasi-compact. $\Box$ \\

In the Appendix we shall give another interesting proof to the quasi-compactness of $\Sp(R)$ without using the flat or Zariski topologies.\\

\begin{lemma} \label{goodlem} Let $J$ be a proper regular ideal of $R$. Then $J$ is a max-regular ideal of $R$ if and only if $R/J$ has no non-trivial idempotents. In particular, if $J$ is a max-regular ideal of $R$ then $\V(J)$ is a connected subset of $\Spec(R)$ w.r.t. the both flat and Zariski topologies.\\
\end{lemma}

{\bf Proof.} If the idempotents of $R/J$ are trivial then it is easy to see that $J$ is a max-regular ideal of $R$. Conversely, suppose $J$ is a max-regular ideal of $R$. Let $\overline{x}=x+J$ be an idempotent element of $R/J$ with $x\in R$. Therefore $x-x^{2}\in J$. We have either $J+(x)=R$ or $J+(x)\neq R$. If $J+(x)=R$ then we may write $1=rx+r_{1}e_{1}+...+r_{m}e_{m}$ with $r,r_{i}\in R$ and $e_{i}\in J$. Thus $1-x=r(x-x^{2})+r'_{1}e_{1}+...+r'_{m}e_{m}$ which belongs to $J$ where $r'_{i}=r_{i}(1-x)$ and so $\overline{x}=1+J$ in this case. But if $J+(x)\neq R$. Since $J$ is a regular ideal thus there are idempotent elements $e'_{1},...,e'_{n}\in J$ and also elements $s_{1},...,s_{n}\in R$ such that $x-x^{2}=s_{1}e'_{1}+...+s_{n}e'_{n}$ we have then $(x-x^{2})\prod\limits_{i=1}^{n}(1-e'_{i})=0$.
By induction on $n$, we observe that there is an idempotent element $e\in J$ such that $\prod\limits_{i=1}^{n}(1-e'_{i})=1-e$. This implies that $(1-e)x$ is an idempotent element. Write $x= ex+(1-e)x$ which belongs to the regular ideal $J'=J+\big((1-e)x\big)$. We have $ J'=J+(x)\neq R$ therefore $J= J'$ since $J$ is a max-regular ideal. Thus $\overline{x}=0$ in this case. The last part of the assertion implies from Corollary \ref{coro 2}. $\Box$ \\

Now we are ready to prove the Theorem:\\

{\bf Proof of Theorem \ref{prop 9}.} Consider the map
$\psi: \Spec(R)\rightarrow \Sp(R)$ which was defined in the proof of Lemma \ref{pierce lemma}. Let $C$ be a connected component of $\Spec(R)$ with respect to the flat topology. Then $\psi(C)$ is a connected subset of $\Sp(R)$. By Lemma \ref{pierce lemma}, $\Sp(R)$ is totally disconnected, so $\psi(C)=\{J\}$ for some max-regular ideal $J$ of $R$. We have then $C\subseteq \psi^{-1}(\{J\})=\V(J)$. By Lemma \ref{goodlem}, $\V(J)$ is a connected subset of $\Spec(R)$. Thus the inclusion $C\subseteq V(J)$ is not strict because $C$ is a connected component. Therefore $C=V(J)$.
Conversely, assume that $J$ is a max-regular ideal of $R$. Again by Lemma \ref{goodlem}, $V(J)$ is a connected subset of $\Spec(R)$ thus it is contained in a flat connected component $C$  of $\Spec(R)$. By the above paragraph, $C= V(J')$ for some max-regular ideal $J'$ of $R$. But $\V(J)\subseteq V(J')$ implies that $\sqrt{J'}\subseteq \sqrt{J}$. We have then $J'\subseteq J$ since $J'$ is a regular ideal. In fact $J'=J$ since $J'$ is a max-regular ideal. Thus $\V(J)=C$ is a flat connected component. The whole of this argument establishes the equivalence $\mathbf{(a)}\Leftrightarrow\mathbf{(c)}$. With a similar argument as above we obtain the equivalence
$\mathbf{(b)}\Leftrightarrow\mathbf{(c)}$. $\Box$ \\

If we equip the set $X=\Spec(R)$ with the flat (resp. Zariski) topology then we shall denote it by $X_{F}$ (resp. $X_{Z}$). By using the flat topology and Hochster's seminal work \cite{Hochster}, a remarkable result in commutative algebra come to light:\\

\begin{theorem}\label{coro 12453} For a given ring $R$ then there exist a ring $A$ and a bijection $X=\Spec(R)\rightarrow Y=\Spec(A)$ which induces the homeomorphisms $X_{F}\rightarrow Y_{Z}$ and $X_{Z}\rightarrow Y_{F}$. In particular, the prime ideals of $A$ have precisely the reverse order of the primes of $R$.\\
\end{theorem}

{\bf Proof.} By Theorem \ref{remark 122} and Lemma \ref{lemma 4465}, $X_{F}$ is a spectral space. Therefore, by \cite[Theorem 6]{Hochster}, there is a ring $A$ and a homeomorphism $\psi:X_{F}\rightarrow Y_{Z}$ where $Y=\Spec(A)$. Let $E$ be a flat closed subset of $Y$. Using Theorem \ref{prop 10}, then $\psi^{-1}(E)$ is stable under specialization. Moreover $\psi^{-1}(E)$ is constructible closed since every spectral map (continuous and quasi-compact) between the spectral spaces is constructible continuous. Therefore, by Theorem \ref{prop 10}, $\psi^{-1}(E)$ is a Zariski closed subset of $X$. By using the same argument for the Zariski closed subsets of $X$ we conclude that $\psi:X_{Z}\rightarrow Y_{F}$ is a homeomorphism. $\Box$ \\

\section{Noetherian aspects of the flat topology}

\begin{lemma}\label{lemma 342} Let $R$ be a ring and let $f\in R$. Then the canonical ring map $\pi: R\rightarrow R_{f}$ is of finite presentation.\\
\end{lemma}

{\bf Proof.} The map $\pi$ induces the following surjective morphism of $R-$algebras $\psi: R[x]\rightarrow R_{f}$ given by $x\rightsquigarrow 1/f$. It suffices to show that $\Ker(\psi)=(fx-1)$. Let $I=(fx-1)$. The image of $f$ under the canonical map $\eta:R\rightarrow R[x]/I$ is invertible since $(f+I)(x+I)=1+I$. Thus, by the universal property of the localizations, there is a (unique) ring map $\lambda: R_{f}\rightarrow R[x]/I$ such that $\eta=\lambda\circ\pi$.  Indeed,
$\lambda(r/f^{n})= rx^{n}+I$.
Therefore the canonical map $\pi': R[x]\rightarrow R[x]/I$ factors as $\pi'=\lambda\circ\psi$. This, in particular, implies that $\Ker(\psi)=I$.  $\Box$ \\

The following result provides various characterizations for the noetherianess of the flat topology.\\

\begin{theorem}\label{theorem 1} The following conditions are equivalent.\\
$\mathbf{(i)}$ The set $\Spec(R)$ equipped with the flat topology is noetherian.\\
$\mathbf{(ii)}$ Every irreducible and closed subset of  $\Spec(R)$ w.r.t. the flat topology is standard Zarsiki open.\\
$\mathbf{(iii)}$ For each prime ideal $\mathfrak{p}$ of $R$ there is some $f\in R\setminus\mathfrak{p}$ such that the canonical map $R_{f}\rightarrow R_{\mathfrak{p}}$ given by $r/f^{n}\rightsquigarrow r/f^{n}$ is bijective.\\
$\mathbf{(iv)}$ For each prime ideal $\mathfrak{p}$ of $R$ the canonical ring map $R\rightarrow R_{\mathfrak{p}}$ is of finite presentation.\\
$\mathbf{(v)}$ The Zariski opens of $\Spec(R)$ are stable under the arbitrary  intersections.\\
$\mathbf{(vi)}$ For every non-empty family $\{\mathfrak{p}_{i}\}$ of prime ideals of $R$ and for each prime ideal $\mathfrak{p}$ of $R$ if
$\bigcap\limits_{i}\mathfrak{p}_{i}\subseteq\mathfrak{p}$ then $\mathfrak{p}_{j}\subseteq\mathfrak{p}$ for some $j$.\\
\end{theorem}

{\bf Proof.} $\mathbf{(i)}\Rightarrow\mathbf{(ii)}:$ Let $Z$ be an irreducible and closed subset of  $\Spec(R)$ with respect to the flat topology. Then $U=\Spec(R)\setminus Z$ is quasi-compact since $\Spec(R)$ is noetherian. Therefore $U=\V(I)$ where $I=(f_{1},...,f_{n})$ is a finitely generated ideal of $R$. Thus $Z=D(f_{i})$ for some $i$ since $Z$ is irreducible. \\
$\mathbf{(ii)}\Rightarrow\mathbf{(i)}:$ It suffices to show that every flat open subset of $\Spec(R)$ is quasi-compact. Let $\mathcal{S}$ be the set of flat open subsets which are not quasi-compact. Let $\mathcal{C}$ be a totally ordered subset in
$\mathcal{S}$ with respect to the inclusion relation and let $W=\bigcup\limits_{U\in\mathcal{C}}U$. The set
$W$ is not quasi-compact. Because if it is quasi-compact then we may find a finite subset $\{U_{1},...,U_{n}\}$ of elements of $\mathcal{C}$ such that $W=\bigcup\limits_{i=1}^{n}U_{i}$. This implies that $W=U_{j}$ for some $j$ since $\mathcal{C}$ is totally ordered. But this is a contradiction since $U_{j}$ is not quasi-compact. Thus $W\in\mathcal{S}$. Therefore every chain has an upper bound in $\mathcal{S}$ and so by the Zorn's Lemma, the family $\mathcal{S}$ has at least a maximal element $U$. Let $E=\Spec(R)\setminus U$. We easily observe that $E$ is flat irreducible. Because suppose $E=E_{1}\cup E_{2}$ where $E_{1}$ and $E_{2}$ are proper flat closed subsets in $E$. By the maximality of $U$, we observe that $U_{i}=\Spec(R)\setminus E_{i}$ are quasi-compact for $i=1,2$. The quasi-compact flat open subsets of $\Spec(R)$  are precisely of the form $\V(I)$ where $I$ is a finitely generated ideal of $R$. Therefore there are finitely generated ideals $I_{i}$ in $R$ such that $U_{i}=\V(I_{i})$. Since $U=U_{1}\cup U_{2}=\V(I_{1}I_{2})$ and $I_{1}I_{2}$ is a finitely generated ideal of $R$. This implies that $U$ is quasi-compact which is a contradiction. Therefore $E$ is flat irreducible. By the hypothesis, there is some element $f\in R$ such that $E=D(f)$, and so $U=\V(f)$. This again means that $U$ is quasi-compact which is a contradiction. Therefore $\mathcal{S}=\emptyset$.\\
$\mathbf{(ii)}\Rightarrow\mathbf{(iii)}:$ Let $\mathfrak{p}$ be a prime ideal of $R$. By the hypothesis, $\Lambda(\mathfrak{p})=D(f)$ for some $f\in R\setminus\mathfrak{p}$. We show that the canonical ring map $\phi:R_{f}\rightarrow R_{\mathfrak{p}}$ given by $r/f^{n}\rightsquigarrow r/f^{n}$ is bijective. The image of each element $s\in R\setminus\mathfrak{p}$ under the canonical ring map $\pi:R\rightarrow R_{f}$ is invertible in $R_{f}$. Because if the ideal generated by $s/1$ is a proper ideal of $R_{f}$ then it is contained in a prime ideal $\mathfrak{q}R_{f}$ of $R_{f}$ where $\mathfrak{q}$ is a prime ideal of $R$ with $\mathfrak{q}\in D(f)$. But this is a contradiction since $s\in\mathfrak{q}\subseteq\mathfrak{p}$. Therefore, by the universal property of the localizations, there is a (unique) ring map $\psi:R_{\mathfrak{p}}\rightarrow R_{f}$ such that $\pi=\psi\circ \pi'$ where $\pi':R\rightarrow R_{\mathfrak{p}}$ is the canonical map. One can easily verify that $\phi\circ\psi$ and $\psi\circ\phi$ are the identity.\\
$\mathbf{(iii)}\Rightarrow\mathbf{(iv)}:$ It is obvious from Lemma \ref{lemma 342}. \\
$\mathbf{(iv)}\Rightarrow\mathbf{(ii)}:$ By Lemma \ref{lemma 4465} and Corollary \ref{coro 1}, a subset of $\Spec(R)$ with respect to the flat topology is  irreducible and closed if and only if it is of the form $\Lambda(\mathfrak{p})=\{\mathfrak{q}\in\Spec(R) : \mathfrak{q}\subseteq\mathfrak{p}\}$ where $\mathfrak{p}$ is a prime ideal of $R$. On the other hand, every flat ring map which is also of finite presentation then it induces a Zariski open map between the corresponding spectra, see \cite[Tag 00I1]{Johan}. It follows that $\Lambda(\mathfrak{p})$
is a Zariski open subset of $\Spec(R)$, and consequently it is standard Zariski open.\\
$\mathbf{(ii)}\Rightarrow\mathbf{(v)}:$ It suffices to show that the intersection of every family of standard Zariski open subsets is Zariski open. Let $\mathfrak{p}\in\bigcap\limits_{\alpha}D(f_{\alpha})$ where for each $\alpha$, $f_{\alpha}\in R$. Therefore  $\Lambda(\mathfrak{p})\subseteq\bigcap\limits_{\alpha}D(f_{\alpha})$. By the hypothesis, $\Lambda(\mathfrak{p})$ is Zarsiki open.\\
$\mathbf{(v)}\Rightarrow\mathbf{(vi)}:$ Suppose for each $i$, there is an element $f_{i}\in\mathfrak{p}_{i}$ such that $f_{i}\notin\mathfrak{p}$. By the hypothesis, $\bigcap\limits_{i}D(f_{i})$ is Zariski open. It is also contains $\mathfrak{p}$. Thus there exists an element $f\in R$ such that $\mathfrak{p}\in D(f)\subseteq\bigcap\limits_{i}D(f_{i})$.
There is some $j$ such that $\mathfrak{p}_{j}\in D(f)$ because  $\bigcap\limits_{i}\mathfrak{p}_{i}\subseteq\mathfrak{p}$. But $D(f)\subseteq D(f_{j})$. This is a contradiction.\\
$\mathbf{(vi)}\Rightarrow\mathbf{(v)}:$ Let $\{E_{\alpha}\}$ be a family of Zarsiki closed subsets of  $\Spec(R)$. For each $\alpha$ there is an ideal $I_{\alpha}$ in $R$ such that $E_{\alpha}=\V(I_{\alpha})$. By the hypothesis,  $\V(J)\subseteq\bigcup\limits_{\alpha}E_{\alpha}$ where $J=\bigcap\limits_{\alpha}\sqrt{I_{\alpha}}$. The reverse inclusion $\bigcup\limits_{\alpha}E_{\alpha}\subseteq\V(J)$
holds in general. Because for each $\alpha$, $J\subseteq\sqrt{I_{\alpha}}$, therefore $E_{\alpha}=\V(I_{\alpha})=\V(\sqrt{I_{\alpha}})\subseteq\V(J)$. \\
$\mathbf{(v)}\Rightarrow\mathbf{(ii)}:$ For each prime ideal $\mathfrak{p}$ of $R$, $\Lambda(\mathfrak{p})=\bigcap\limits_{f\in R\setminus\mathfrak{p}}D(f)$. Therefore, by the hypothesis, $\Lambda(\mathfrak{p})$ is Zariski open and so it is standard Zariski open.  $\Box$ \\

\begin{lemma}\label{lemma 4562} Let $S$ be a multiplicative subset of a ring $R$ such that the canonical ring map $R\rightarrow S^{-1}R$ is injective and of finite type. Then $S^{-1}R$ is canonically isomorphic to $R[x]/(fx-1)$ for some $f\in R$.\\
\end{lemma}

{\bf Proof.} There exists a finite subset $\{r_{1}/f_{1},...,r_{n}/f_{n}\}$ of $S^{-1}R$ such that each element $r/s\in S^{-1}R$ can be written as $r/s=g(r_{1}/f_{1},...,r_{n}/f_{n})$ for some $g(x_{1},...,x_{n})\in R[x_{1},...,x_{n}]$. It follows that there is a natural number $N$ such that $r/s=r'/f^{N}$ where $f=f_{1}f_{2}...f_{n}$. This means that the canonical map $\phi: R_{f}\rightarrow S^{-1}R$ given by $r/f^{d}\rightsquigarrow r/f^{d}$ is surjective. By the hypothesis it is also injective. Thus the assertion implies from the proof of Lemma \ref{lemma 342}. $\Box$ \\

\begin{corollary} Let $R$ be a domain. Then $\Spec(R)$ equipped with the flat topology is noetherian if and only if for each prime ideal $\mathfrak{p}$ of $R$ the canonical map $R\rightarrow R_{\mathfrak{p}}$ is of finite type.\\
\end{corollary}

{\bf Proof.}  The implication ``$\Rightarrow$" implies from Theorem \ref{theorem 1}.  For the converse apply Lemma \ref{lemma 4562} and Theorem \ref{theorem 1}. $\Box$ \\

The noetherianess of the infinite spectra has a wild nature in the passage between the flat and Zariski topologies. More precisely:\\

\begin{theorem}\label{theorem 4537} Let $R$ be a ring such that $\Spec(R)$ is noetherian w.r.t. the both flat and Zariski topologies. Then $R$ has a finitely many prime ideals.\\
\end{theorem}

{\bf Proof.} Let $\mathfrak{p}$ be a prime ideal of $R$. By Theorem \ref{theorem 1}, there is an element $f\in R\setminus\mathfrak{p}$ such that $\Lambda(\mathfrak{p})=D(f)$. On the other hand, $\Spec(R)\setminus\V(\mathfrak{p})$ is quasi-compact since $\Spec(R)$ is  noetherian with respect to the Zariski topology. Thus there is a finitely generated ideal $I$ of $R$ such that $\V(\mathfrak{p})=\V(I)$.
It follows that $\{\mathfrak{p}\}=D(f)\cap\V(I)$. Therefore, by Remark \ref{remark 432}, $\Spec(R)$ is discrete w.r.t. the constructible topology. The underlying set of a quasi-compact discrete space is finite, hence $\Spec(R)$ is a finite set.  $\Box$ \\

\begin{proposition}\label{coro 13456} Let $R$ be a ring such that $X=\Spec(R)$ is noetherian w.r.t. the flat topology. Then $R$ satisfies the descending chain condition on prime ideals.\\
\end{proposition}

{\bf Proof.} It is an immediate consequence of Theorem \ref{coro 12453}. $\Box$ \\

\begin{remark} Theorem \ref{coro 12453}, in particular, provides many examples of rings with infinitely many prime ideals whose spectra equipped with the flat topology are noetherian. Specially, for every notherian ring $R$, then there is a ring $A$ such that the prime spectrum of $A$ equipped with the flat topology is noetherian and the set $\Spec(A)$ is in bijection with $\Spec(R)$. \\
\end{remark}

\section{Appendix}

Here we give a second proof to the quasi-compactness of the pierce spectrum (Lemma \ref{pierce lemma}). Let $\mathscr{C}=\{U_{e_{i}} : i\in I\}$ be an open covering for $\Sp(R)$ consisting of basis opens which does not have any finite refinement. Consider the binary operation $\wedge$ over $I(R)$ which is defined as $e\wedge e'=e+e'-ee'$.  One can easily verify that this operation is associative and $U_{e}\cup U_{e'}=U_{e\wedge e'}$. The regular ideal $K$ of $R$ generated by $\{e_{i} : i\in I\}$ is proper. Since otherwise we may write $1=r_{1}e_{i_{1}}+...+r_{n}e_{i_{n}}$. This implies that $(1-e_{i_{1}})...(1-e_{i_{n}})=0$. Therefore $e_{i_{1}}\wedge...\wedge e_{i_{n}}=1$. This means that $\Sp(R)=U_{e_{i_{1}}}\cup...\cup U_{e_{i_{n}}}$.
But this is a contradiction since the covering  $\mathscr{C}$ does not have any finite refinement. Thus $K$ is contained in a max-regular ideal $J$. Therefore $e_{i}\in J$ for all $i\in I$, but this is also a contradiction. $\Box$ \\

\textbf{Acknowledgements.} The author would like to give sincere thanks to Professor Marco Fontana for many valuable discussions which the author had with him during the writing the present article and also for his informing about their earlier work. \\


\begin{thebibliography}{10}
\bibitem{Johan}
de Jong A.J., et al., The Stacks Project, see http://stacks.math.columbia.edu., (2019).
\bibitem{Fontana}
Dobbs D., Fontana M., and Papick I. On the flat spectral topology, Rendiconti di Mathematica. (4) 1981, Vol. 1, Serie VII.
\bibitem{Hochster}
Hochster, M. Prime ideal structure in commutative rings, Trans. Amer.
Math. Soc. 142 (1969), 43-60.
\bibitem{Ebulfez 25}
Tarizadeh, A. On the separability of the underlying space of a scheme, submitted,  arXiv:1803.04817v1 [math.AC].
\bibitem{Tarizadeh 2020}
Tarizadeh, A. A note on finitely generated flat modules, submitted, arXiv:1804.03007v1 [math.AC].
\end{thebibliography}
\end{document}